\documentclass[12pt]{amsart}
\usepackage{amssymb,latexsym}
\newdimen\AAdi%
\newbox\AAbo%
\def\AAk#1#2{\s_etbox\AAbo=\hbox{#2}\AAdi=\wd\AAbo\kern#1\AAdi{}}%
\def\AAr#1#2#3{\s_etbox\AAbo=\hbox{#2}\AAdi=\ht\AAbo\raise#1\AAdi\hbox{#3}}%
\font\tenmsb=msbm10 at 12pt \font\sevenmsb=msbm7 at 8pt
\font\fivemsb=msbm5 at 6pt
\newfam\msbfam
\textfont\msbfam=\tenmsb \scriptfont\msbfam=\sevenmsb
\scriptscriptfont\msbfam=\fivemsb
\def\Bbb#1{{\tenmsb\fam\msbfam#1}}

\begin{document}

\newcommand{\rD}{D \hskip -2.8mm \slash}
\newcommand{\D}{\partial \hskip -2.2mm \slash}
\newcommand{\e}{\overline{\varepsilon}}
\renewcommand{\theequation}{\thesection.\arabic{equation}}
\newcommand{\wb}{\widetilde{\nabla}_{e_\beta}}
\newcommand{\wa}{\widetilde{\nabla}_{e_\alpha}}

\newtheorem{thm}{Theorem}
\newtheorem{lem}{Lemma}
\newtheorem{cor}{Corollary}
\newtheorem{rem}{Remark}
\newtheorem{pro}{Proposition}
\newtheorem{defi}{Definition}
\newcommand{\noi}{\noindent}
\newcommand{\dis}{\displaystyle}
\newcommand{\mint}{-\!\!\!\!\!\!\int}
\newcommand{\ba}{\begin{array}}
\newcommand{\ea}{\end{array}}
\def \tf{\tilde{f}}
\def\cqfd{%
\mbox{ }%
\nolinebreak%
\hfill%
\rule{2mm} {2mm}%
\medbreak%
\par%
}
\def \pr {\noindent {\it Proof.} }
\def \rmk {\noindent {\it Remark} }
\def \esp {\hspace{4mm}}
\def \dsp {\hspace{2mm}}
\def \ssp {\hspace{1mm}}
\def\n{\nabla}
\def\RR{\Bbb R}\def\R{\Bbb R}
\def\C{\Bbb C}
\def\B{\Bbb B}
\def\N{\Bbb N}
\def\Q{\Bbb Q}
\def\Z{\Bbb Z}
\def\EE{\Bbb E}
\def\H{\Bbb H}
\def\SS{\Bbb S}\def\S{\Bbb S}
\def \c {{\bf C}}
\def \Q {{\bf Q}}
\def \a {\alpha}
\def \b {\beta}
\def \d {\delta}
\def \e {\epsilon}
\def \G {\Gamma}
\def \g {\gamma}
\def \l {\lambda}
\def \L {\Lambda}
\def \O {\Omega}
\def \om {\omega}
\def \o{\omega}
\def \s {\sigma}
\def \S {\Sigma}
\def \t {\theta}
\def \z {\zeta}
\def \vp {\varphi}
\def \vt {\vartheta}
\def \ve {\varepsilon}
\def \i {\infty}
\def \ds {\displaystyle}
\def \oo {\overline{\Omega}}
\def \ov {\overline}
\def \bd {\bigtriangledown}
\def \U {\bigcup}
\def \un {\underline}
\def \h {\hspace{.5cm}}
\def \hs {\hspace{2.5cm}}
\def \v {\vspace{.5cm}}
\def \mi {M_{i}}
\def \ra {\longrightarrow}
\def \Ra {\Longrightarrow}
\def \rw {\rightarrow}
\def \bs {\backslash}
\def \rn {{\bf R}^n}
\def \h* {\hspace*{1cm}}
\def\la{\langle}
\def\ra{\rangle}
\def\cal{\mathcal}

\title[ Dirac equations ]{Nonlinear Dirac equations on Riemann surfaces}

\author{
Qun Chen, J\"urgen Jost and Guofang Wang }

\thanks{The research of QC is partially supported by  NSFC (Grant No.10571068) and SRF for
ROCS, SEM; he also thanks the Max Planck Institute for Mathematics
in the Sciences for good working conditions during his visit.}

\address{School of Mathematics and Statistics\\ Wuhan University\\Wuhan 430072, China } \email{qunchen@whu.edu.cn}

\address{Max Planck Institute for Mathematics in the Sciences\\Inselstr. 22\\D-04103 Leipzig, Germany} \email{jjost@mis.mpg.de}

\address{Faculty of Mathematics, University Magdeburg, D-39016, Magdebrug, Germany}
\email{gwang@math.uni-magdeburg.de}


\begin{abstract}
We develop analytical methods for nonlinear Dirac equations. Examples
of such equations include
Dirac-harmonic maps with curvature term and the equations describing  the generalized Weierstrass representation of surfaces in
three-manifolds. We provide the key analytical steps, i.e., small
energy regularity and
removable singularity theorems and energy identities for solutions.
\end{abstract}
\date{April 2, 2007}
\maketitle

{\small
 } \vskip24pt

\section{Introduction}
Dirac type equations on Riemann surfaces are ubiquitous in geometry
as they constitute the most basic first order system of elliptic
equations. The first examples are of course the Cauchy-Riemann
equations. These are linear, but other examples are typically of the
nonlinear type
\begin{equation}
\label{1.0} \D\psi=H_{jkl}\la\psi^j,\psi^k\ra\psi^l,
\end{equation}
with a notation to be explained shortly. The linear operator on the
left hand side is, of course, the Dirac operator whereas the
nonlinearity on the right hand side is  cubic. As we shall see,
this type of nonlinearity on one hand arises naturally in geometry,
because (\ref{1.0}) is conformally invariant
and on the other hand, from the analytical side, it presents a
borderline case where standard linear methods fail to apply (again,
because it is conformally invariant), but an
analytical treatment nevertheless is still possible by utilizing the
structure of the equation in a more sophisticated manner. That is,
analytical methods need to be supplemented by geometric insights. This
frame makes (\ref{1.0}) particularly attractive. In the present paper,
we develop a systematic and general treatment of the key steps of the
nonlinear analysis.\\

Let us now describe the underlying geometric structure in more detail.
Let $(M,g)$ be a Riemann surface with a fixed spin structure,
and denote  the spin bundle by $\S$. On $\Sigma$, there is a Hermitian metric
$\la\cdot,\cdot\ra$ compatible with the spin connection $\nabla$ on
$\S$. For any orthonormal basis $\{e_\alpha, \alpha=1,2\}$ on $M$,
the (Atiyah-Singer) Dirac operator is defined by
$\D:=e_\alpha\cdot\nabla_{e_\alpha}$,  where $\cdot$ stands for the
Clifford multiplication. In this paper, we use the summation
convention.

We consider (\ref{1.0}), that is,
\begin{equation}
\label{1.1} \D\psi=H_{jkl}\la\psi^j,\psi^k\ra\psi^l,
\end{equation}
where $\psi=(\psi^1,\psi^2,\cdots,\psi^n)\in\Gamma (\S^n),$
$\S^n:=\overbrace{\S\times\cdots\times\S}\limits^n$, $n\in \Z_+$ and
$H_{jkl}=(H^1_{jkl},H^2_{jkl},$  $\cdots,H^n_{jkl})\in C^1(M,\R^n).$
   Denote   $h_0:=\max\{|H^i_{jkl}|(x);$ $ x\in M,
   i,j,k,l=1,2,\cdots,n\}$,
$h_1:=\max\{|\n H^i_{jkl}|(x); x\in M, i,j,k,l=1,2,\cdots,n\}$ and
$|\psi|:=(\sum\limits_{i=1}^n\la\psi^i,\psi^i\ra)^{1/2}$. We note
that (\ref{1.1})  is conformally invariant.

Let us now discuss examples where (\ref{1.1}) arises. In fact, we
have been led to it through our study of  Dirac-harmonic maps with
curvature term (c.f. \cite{CJW2}) which in turn were derived from
the nonlinear supersymmetric $\sigma$-model of quantum field theory
where Dirac type equations describe fermionic particles.  Let $\phi$
be a smooth map from $M$ to a Riemannian manifold $(N,h)$ of
dimension $n \ge 2$ and $\phi^{-1} TN$ the pull-back bundle of $TN$
by $\phi$. On the twisted bundle $\Sigma \otimes \phi^{-1} TN$,
there is a metric (also denoted by $\la \cdot,\cdot\ra$) and a
natural connection $\widetilde{\n}$ induced from  those on $\Sigma $
and $\phi^{-1} TN$.  In local coordinates $\{x_\alpha\}$ and
$\{y^i\}$ on $M$ and $N$ respectively, a section $\psi$ of $\Sigma
\otimes \phi^{-1} TN$ takes the form
\[\psi(x)= \psi^j (x)\otimes\partial_{y^j}(\phi(x))\]
and $\widetilde{\n}$ can be written as
\[\widetilde{\nabla} \psi(x)= \n \psi^i(x)\otimes \partial_{y^i}(\phi(x))
+ \Gamma^i_{jk}(\phi(x)) \n \phi^j(x) \psi^k(x)\otimes\partial_{y^i}
(\phi(x)),\] where $\psi^i\in \Gamma(\Sigma)$, $\{\partial_{y^j}\}$
is the natural local basis on $N$ and $\{\Gamma^i_{jk}\}$ stands for
the Christoffel symbols of $N$. The {\it Dirac operator along the
map $\phi$} is defined as
 \begin{eqnarray*}
\rD\psi&:=& e_\a\cdot \widetilde{\nabla}_{e_\a}
\psi\nonumber \\
&=&
\partial \hskip -2.2mm \slash
 \psi^i(x)\otimes\partial_{y^i}(\phi(x))
+ \Gamma^i_{jk}(\phi(x)) \n_{e_\a} \phi^j(x) e_\a\cdot
\psi^k(x)\otimes\partial_{y^i}(\phi(x)).
\end{eqnarray*}
In \cite{CJW2}, we considered the following functional:
\begin{equation}
\label{b4} L_c(\phi,\psi):=\frac 1 2
\int_M[|d\phi|^2+\la\psi,\rD\psi\ra-\frac 1 6 R_{ikjl}\la
\psi^i,\psi^j \ra \la \psi^k,\psi^l \ra].
\end{equation}
We call critical points $(\phi,\psi)$ of $L_c$  {\it
Dirac-harmonic maps with curvature term}. This functional is dictated
by the supersymmetry requirements of the $\sigma$-model in superstring theory. The
difference is that here the components of $\psi$ are ordinary spinor
fields on $M$, while in physics they take values in a Grassmann
algebra. The Euler-Lagrange equations of the functional $L_c$ are
 (see \cite{CJW2} for details):
\begin{equation}
\label{3.2.11}
\tau^i(\phi)-\frac{1}{2}R^i\hskip0.000001mm_{lmj}\la\psi^m,\nabla\phi^l\cdot\psi^j\ra+\frac{1}{12}h^{ip}R_{mkjl;p}
\la\psi^m,\psi^j\ra\la\psi^k,\psi^l\ra=0,
\end{equation}
\begin{equation}
\label{3.2.1}
\rD\psi^i=-\frac{1}{3}R^i\hskip0.000001mm_{jkl}\la\psi^j,\psi^k\ra\psi^l,
\qquad i=1,2,\cdots,n,
\end{equation}
where $\tau(\phi)$ is the tension field of $\phi$,
$R^i\hskip0.000001mm_{jkl}$ stands for a component of the curvature
tensor of $N$ and $R_{mkjl;p}$ denotes the covariant derivative of
$h_{mi}R^i\hskip0.000001mm_{kjl}$ with respect to $\partial_{y^p}$.

In particular, if $\phi$ is a constant map, then (\ref{3.2.1})
becomes
\begin{equation}
\label{3.2.2}
\D\psi^i=-\frac{1}{3}R^i\hskip0.000001mm_{jkl}\la\psi^j,\psi^k\ra\psi^l,\qquad
i=1,2,\cdots,n,
\end{equation}
which is a Dirac equation of type (\ref{1.1}).

Another more classical example of an equation of type  (\ref{1.1})
comes from the geometry of surfaces in three-manifolds  via the
generalized Weierstrass representation as we shall now explain. When
$n=1$, it takes the form
\begin{equation}
\label{3.2.3} \D\psi=H|\psi|^2\psi \end{equation}
 for $H\in
C^1(M)$, $\psi\in \Gamma(\Sigma)$. A similar equation was considered by Ammann and Humbert in \cite{AH}
when they studied  the first conformal Dirac eigenvalue. See also \cite{A} for a Yamabe type problem.

Recall that the classical
Weierstrass formula represents minimal surfaces $X$ immersed in
$\R^3$ in terms of a holomorphic 1-form and a meromorphic function.
The generalized Weierstrass representation was found to express a
general surface immersed in $\R^3$ (as well as in $\R^4$ and some
three-dimensional Lie groups) by Dirac equations (see e.g.
\cite{Fr}, \cite{Ke}, \cite{Ta2}). For the ambient space $\R^3$, a
surface $X:M\to \R^3$ is represented by
$$X=Re \int(i(\psi^2_1+\bar{\psi}^2_2), \bar{\psi}^2_2-\psi^2_1, 2\psi_2\bar{\psi}_2),$$
where $\psi:=\left(\begin{matrix} \ds \psi_1
\\ \ds \psi_2\end{matrix}\right)
:\R^2\to \C^2$ satisfies the following equation:
\begin{equation}
\label{1.1.2}
\left[2
\ds \left(\begin{matrix}0& \bar\partial\\ -\partial &0
\end{matrix}\right)+\ds \left(\begin{matrix}U& 0\\ 0 &U
\end{matrix}\right)\right]
\left(\begin{matrix} \ds
 \psi_1\\ \ds \psi_2
\end{matrix}\right)=0,
\end{equation}
with $U=H|\psi|^2.$  On the Euclidean plane $\R^2$ the spin
structure is unique, and the spinor bundle $\S$ is trivial. By
choosing a representation of $\{e_\alpha\}$ as
\[\sigma_1=\left(\begin{matrix}0& 1\\ -1&0 \end{matrix}\right),
\quad \sigma_2=\left(\begin{matrix}0& i\\ i&0
\end{matrix}\right),\qquad i:=\sqrt{-1}\] the Dirac operator can be expressed as
\begin{eqnarray}
\label{1.4}
\partial \hskip
-2.2mm \slash\psi&=&\ds \left(\begin{matrix}0& 1\\ -1&0
\end{matrix}\right) \left(\begin{matrix} \ds \partial_x
 \psi_1\\ \ds \partial_x  \psi_2
\end{matrix}\right)+
\left(\begin{matrix}0& i\\ i&0 \end{matrix}\right)
\left(\begin{matrix} \ds \partial_y  \psi_1 \\
\ds\partial_y  \psi_2
\end{matrix}\right)\nonumber\\
&=&2 \ds \left(\begin{matrix}0& \bar\partial\\ -\partial &0
\end{matrix}\right)
\left(\begin{matrix} \ds
 \psi_1\\ \ds \psi_2
\end{matrix}\right),
\end{eqnarray}
where $\psi:=\left(\begin{matrix} \ds \psi_1
\\ \ds \psi_2\end{matrix}\right)
:\R^2\to \C^2$ is a spinor field, $\partial:=\frac 1 2
(\partial_x-i\partial_y)$, $\bar \partial:=\frac 1 2
(\partial_x+i\partial_y)$. Therefore, the equation (\ref{1.1.2}) can
be written as (\ref{3.2.3}) which is a special case of (\ref{1.1})
with $n=1$. Similar types of equations will be discussed in section
5.


We now turn to the analytical aspects and introduce the following
energy functional:
\begin{equation}
\label{1.2} E(\psi):=\int_M |\psi|^4.
\end{equation}
For the analysis of the equation (\ref{1.1}), we use
$\|\psi\|_{L^4}$, instead of $\|\nabla\psi\|_{L^{\frac 4 3}}$, as
the energy functional of $\psi$, since the former is strong enough
to get various estimates in most cases, as one can see in
\cite{CJLW1} and \cite{CJLW2}.

For a solution $\psi$ of (\ref{1.1}), if $\psi\in L^r$ for some
$r>4$, then the standard argument of elliptic regularity theory
implies the smoothness of $\psi$. Furthermore, under the condition
of uniformly bounded $L^r(r>4)$ norms, one has compactness for these
solutions. However, if we only assume the boundedness of the $L^4$
norm of $\psi$, then the compactness is no longer true. One then
naturally considers the blow up phenomenon for the solutions. In
particular, it is interesting to know whether the energy identity
and removable singularity results hold.  $\psi\in L^4$ is the
borderline case for the geometric analysis of the solutions of
(\ref{1.1}), since the standard bootstrap method in elliptic
estimates fails in the first step. It turns out that this can be
overcome by establishing some $L^p$ elliptic boundary estimates and
combining various estimates in delicate ways.

In this paper, our aim is to study properties of the solution space
of this and similar types of equations. We will prove small energy
regularity theorems in which basic apriori estimates for smooth
solutions of (\ref{1.1}) and related types of equations are given.
Then we prove a removable singularity theorem, which provides a key
tool for the blow up analysis of the solutions. Based on these, we
consider the blow up process of solutions under a uniform bound for
the
$L^4$ norm, and we will establish an energy identity for these
solutions. For harmonic maps, holomorphic curves and also maps with
uniformly $L^2$-bounded tension fields, these results are derived in
\cite{SU},  \cite{J2}, \cite{DT}, \cite{Ye}, \cite{Pk} and
\cite{PW}.


We would like to remark that these considerations are closely
related to  the regularity of weak solutions, corresponding to the
well-known case of harmonic maps in dimension two. Based on the
discussion in this paper, we believe that, using an $L^p$ theory for
boundary value problems of the Dirac equations (see \cite{BP} for
the $L^2$ theory), regularity results for weak solutions are
available.

\vskip36pt
\section{Small energy regularity theorem}
\addtocounter{equation}{-9}

 In this section, we will prove a small energy regularity theorem. Since the
problem is local and the equation (\ref{1.1}) is conformally
invariant, without loss of generality, we may assume $M$ to be the
unit disk
$$B=\{(x,y)\in\R^2|x^2+y^2<1\}$$
equipped with the Euclidean metric.

\vskip12pt \noindent {\bf Theorem 2.1 (Small energy regularity
theorem).} {\it There exists a small constant $\varepsilon>0$ such
that for any smooth solution $\psi$ of (\ref{1.1}) satisfying
\begin{equation}
\label{2.1} E(\psi;B):=\int\limits_B|\psi|^4<\varepsilon,
\end{equation}
we have
\begin{equation}
\label{2.2}
\|\psi\|_{B',k,p}\leq C\|\psi\|_{B,0,4},
\end{equation}
$\forall B'\subset \subset B,$ $1<p$ and $k\in\Z_+,$ where
$C=C(B',k,p)>0$ is a constant, and $\|\cdot\|_{B,k,p}$ denotes the
norm in $W^{k,p}(B,\S^n)$.} \vskip12pt

For proving this result, we need the following  $L^p$ boundary
estimates for Dirac operators. This is essentially Lemma 4.8 in
\cite{CJLW1}, but we will give another proof here. \vskip12pt
\noindent {\bf Lemma 2.2.} {\it Suppose $\psi$ is a solution of
\begin{equation}
\label{2.3}
\D\psi=f
\end{equation}
on $B$, with $\psi|_{\partial B}=\varphi$, and $f\in L^p(B,\S^n)$,
$\varphi \in W^{1,p}(\partial B,\S^n)$ for some $p>1$, then
\begin{equation}
\label{2.4}
\|\psi\|_{B,1,p}\leq C(\|f\|_{B,0,p}+\|\varphi\|_{\partial B,1,p}),
\end{equation}
where $C=C(p)>0$ is a constant. }  \vskip12pt \noindent {\it Proof
of Lemma 2.2.}  First, let
$$w(x):=\int\limits_B\Gamma(x-y)f(y)dy$$
be the Dirac-Newton potential of $f$, where
$$\Gamma(x):=-\frac{1}{2\pi} \frac{x}{|x|^2}\hskip0.6mm\cdot$$ is the Green function of $\D$. Using the relation between
 $\Gamma$ and the Green function $\widetilde{\Gamma}$ of the Laplace operator: $\Gamma=\D\widetilde{\Gamma}$,
 one sees that $w=\D\widetilde{w}$, where
$$\tilde{w}(x):=\int\limits_B\widetilde{\Gamma}(x-y)f(y)dy.$$
From the Calderon-Zygmund inequality (see e.g. \cite{GT} Theorem
9.9):
$$\|\nabla^2\tilde{w}\|_{B,0,p}\leq C\|f\|_{B,0,p}.$$
Noting that
$$\nabla_\alpha w=\nabla_\alpha\D\tilde{w}=\sigma_\beta\nabla_\alpha\nabla_\beta\tilde{w},$$
it follows that
\begin{equation}
\label{2.5}
\|\nabla w\|_{B,0,p}\leq C\|f\|_{B,0,p}.
\end{equation}
Second, for any $\xi\in W^{1,p}_0(B,\S^n)$ and $F\in L^p(B,\S^n)$
satisfying
\begin{equation}
\label{2.6}
\D\xi=F,
\end{equation}
we have
\begin{equation}
\label{2.7}
\|\nabla\xi\|_{B,0,p}\leq\|F\|_{B,0,p}.
\end{equation}
In fact, there exists a sequence of $\xi_k\in C^1_0(B,\S^n)$ such
that $\xi_k\to\xi$ in $W^{1,p}$ which implies $\D\xi_k\to \D\xi$ in
$L^p$. Denote $F_k:=\D\xi_k$, since $\xi_k$ has compact support , we
know (see e.g. \cite{MS}) that $\xi_k$ is the Dirac-Newton potential
of $F_k$, and by (\ref{2.5}) we have
$$\|\nabla\xi_k\|_{B,0,p}\leq\|F_k\|_{B,0,p},$$
letting $k\to+\infty$ yields (\ref{2.7}).
\par
Now we extend $\varphi$ to $\tilde{\varphi}$ on $B\setminus B_\delta$ $(0<\delta<\frac 1 2)$ by
$$\tilde{\varphi}(r,\theta):=\varphi(\theta),\quad \delta\leq r\leq 1,\quad \theta\in \partial B.$$
Choose a cut-off function $\eta$ such that $0\leq \eta\leq 1$,
$$
            \eta=\left\{
                \begin{array}{l}
  1\quad r\geq \frac 3 4, \\
  0\quad r\leq\frac 1 2,
                \end{array}
                \right.
  $$
and $|\eta'|\leq 2$, define
$$\hat{\varphi}:=\eta\tilde{\varphi},$$
then $\hat{\varphi}\in W^{1,p}(B,\S^n)$ and $\psi -\hat{\varphi}\in
W_0^{1,p}(B,\S^n).$
\par
From (\ref{2.7}), we have
\begin{eqnarray*}
\|\nabla(\psi-\hat{\varphi})\|_{B,0,p}&\leq & C\|\D(\psi-\hat{\varphi})\|_{B,0,p}\\
&\leq &C(\|f\|_{B,0,p}+\|\D\hat{\varphi}\|_{B,0,p}),
\end{eqnarray*}
which implies
\begin{equation}
\label{2.8}
\|\nabla\psi\|_{B,0,p}\leq C(\|f\|_{B,0,p}+\|\D\hat{\varphi}\|_{B,0,p}+\|\nabla\hat{\varphi}\|_{B,0,p}).
\end{equation}
Note that

\begin{eqnarray*}
\|\nabla\hat{\varphi}\|_{B,0,p}&=& \|\nabla(\eta\tilde{\varphi})\|_{B,0,p}\\
&\leq & C(\|\tilde{\varphi}\|_{B_{\frac 3 4}\setminus B_{\frac 1 2},0,p}+\|\nabla\tilde{\varphi}\|_{B\setminus
B_{\frac 1 2},0,p})\\
&\leq & C\|\varphi\|_{\partial B,0,p}+[\int\limits_{B\setminus B_{\frac 1 2}}(\frac 1 r|\nabla\varphi|)^p]^{\frac 1 p}
\\
&\leq & C(\|\varphi\|_{\partial B,0,p}+\|\nabla\varphi\|_{\partial B,0,p}),
\end{eqnarray*}
namely,
\begin{equation}
\label{2.9} \|\nabla\hat{\varphi}\|_{B,0,p}\leq
C\|\varphi\|_{\partial B,1,p}.
\end{equation}
Similarly,
\begin{equation}
\label{2.10} \|\D\hat{\varphi}\|_{B,0,p}\leq C\|\varphi\|_{\partial
B,1,p}.
\end{equation}
Subsituting (\ref{2.9}) and (\ref{2.10}) into (\ref{2.8}) then yields:
\begin{equation}
\label{2.11}
\|\nabla\psi\|_{B,0,p}\leq C(\|f\|_{B,0,p}+\|\varphi\|_{\partial B,1,p}).
\end{equation}

\noindent By the Poincar\'e inequality,
$$\|\psi-\hat{\varphi}\|_{B,0,p}\leq C \|\nabla(\psi-\hat{\varphi})\|_{B,0,p},$$
hence
\begin{equation}
\label{2.12}
\|\psi\|_{B,0,p}\leq C \|\nabla(\psi-\hat{\varphi})\|_{B,0,p}+\|\hat{\varphi}\|_{B,0,p},
\end{equation}
but $$\|\hat{\varphi}\|_{B,0,p}=\|\eta\tilde{\varphi}\|_{B,0,p}\leq C\|\varphi\|_{\partial B,0,p},$$
putting this and (\ref{2.9}), (\ref{2.11}) into (\ref{2.12}), we finally obtain (\ref{2.4}).
\hfill Q.E.D.

\vskip12pt Now we can give the \vskip12pt \noindent {\it Proof of
Theorem 2.1.}  We first derive the following estimate:
\begin{equation}
\label{2.13}
\|\psi\|_{B',0,p}\leq C \|\psi\|_{B,0,4},
\end{equation}
$\forall B'\subset \subset B$, where $C=C(B',p)>0$ is a constant.
\par For this, we choose a cut-off function $\eta$ such that $0\leq \eta\leq 1$, $\eta|_{B'}\equiv 1$, and
${\rm supp}\eta\subset B$. Denote $\xi:=\eta\psi$, then
\begin{eqnarray}
\label{2.14}
\D\xi&=&\D(\eta\psi)\nonumber\\
&=&\eta\D\psi+\nabla\eta\cdot\psi\nonumber\\
&=&\eta H_{jkl}\la\psi^j,\psi^k\ra\psi^l+\nabla\eta\cdot\psi.
\end{eqnarray}
From Lemma 2.2, for any $1<q<2$,
\begin{eqnarray}
\label{2.15}
\|\xi\|_{B,1,q}&\leq &C\|\eta H_{jkl}\la\psi^j,\psi^k\ra\psi^l+\nabla\eta\cdot\psi\|_{B,0,q}\nonumber\\
&\leq &C(h_0\|\eta|\psi|^3\|_{B,0,q}+\|\psi\|_{B,0,q}).
\end{eqnarray}
Now observing that
\begin{eqnarray}
\label{2.16}
\|\eta|\psi|^3\|_{B,0,q}&=&[\int\limits_B(|\psi|^2|\eta\psi|)^q]^{\frac 1 q}\nonumber\\
&=&(\int\limits_B(|\psi|^{2q}|\xi|^q)^{\frac 1 q}\nonumber\\
&\leq&(\int\limits_B|\psi|^4)^{\frac 1 2}(\int\limits_B|\xi|^{q^*})^{\frac{1}{q^*}},
\end{eqnarray}
where $q^*:=\frac{2q}{2-q}$, putting (\ref{2.16}) into (\ref{2.15}), and using the Sobolev embedding, we have
$$\|\xi\|_{B,0,q^*}\leq C(h_0\|\psi\|^2_{B,0,4}\|\xi\|_{B,0,q^*}+\|\psi\|_{B,0,4}).$$
Thus, if $\varepsilon>0$ is small enough such that
$Ch_0\sqrt{\varepsilon}<\frac 1 2$, then for $\psi$ with
$$\int\limits_B|\psi|^4<\varepsilon$$
we have
$$\|\xi\|_{B,0,q^*}\leq C\|\psi\|_{B,0,4}.$$
Clearly, for any $p>1$, one can find some $q<2$ such that $p=q^*$. This establishes (\ref{2.13}).
\par
Next, since
\begin{eqnarray*}
\int\limits_B |\nabla\xi|^2&=&\int\limits_B|\D\xi|^2\\
&=&\int\limits_B|\eta\D\psi+\nabla\eta\cdot\psi|^2\\
&\leq& C(\int\limits_B|\psi|^6+\int\limits_B|\psi|^2),
\end{eqnarray*}
we have
\begin{eqnarray*}
\|\nabla\xi\|_{B,0,2}&\leq& C(\|\psi\|^3_{B,0,6}+\|\psi\|_{B,0,2})\\
&\leq& C(\|\psi\|^3_{B,0,4}+\|\psi\|_{B,0,4})\\
&\leq& C\|\psi\|_{B,0,4}(1+\|\psi\|^2_{B,0,4})\\
&\leq& C\|\psi\|_{B,0,4},
\end{eqnarray*}
where in the second step we have used (\ref{2.13}) and in the last step we have used (\ref{2.1}).
We then have,
\begin{equation}
\label{2.17}
\|\nabla\psi\|_{B',0,2}\leq C\|\psi\|_{B,0,4},
\end{equation}
where $C>0$ is constant depending only on $h_0$ and $B'$.
Using the Weitzenb\"{o}ck formula and noting that the scalar curvature vanishes in this case, we have
\begin{eqnarray*}
\nabla_\alpha\nabla_\alpha\psi&=&-\D^2\psi\\
&=&-\D(H_{jkl}\la\psi^j,\psi^k\ra\psi^l)\\
&=&-\la\psi^j,\psi^k\ra(\n H_{jkl}\cdot\psi^l)-H_{jkl}(\la\n_{e_\alpha}\psi^j,\psi^k\ra+\la\psi^j,\n_{e_\alpha}\psi^k\ra)(e_\alpha\cdot\psi^l)\\
&&-H_{jkl}H^l_{pqr}\la\psi^j,\psi^k\ra\la\psi^p,\psi^q\ra\psi^r.
\end{eqnarray*}
Therefore, for any $\eta\in C^\infty (B)$,
$$|\Delta(\eta\psi)|\leq C(|\psi|+|\nabla\psi|+|\nabla\psi||\psi|^2+|\psi|^3+|\psi|^5),$$
where $C>0$ is a constant depending only on $\eta, h_0$ and $h_1$, from which we have
\begin{equation}
\label{2.18} \|\eta\psi\|_{2,p}\leq
C(\|\psi\|_{0,p}+\|\nabla\psi\|_{0,p}+\||\nabla\psi||\psi|^2\|_{0,p}
+\||\psi|^3\|_{0,p} +\||\psi|^5\|_{0,p}).
\end{equation}
Using the above estimates (\ref{2.13}) and (\ref{2.17}), we have
\begin{eqnarray*}
\||\nabla\psi||\psi|^2\|_{B',0,\frac 4 3}&\leq&\|\nabla\psi\|_{B',0,2}\|\psi\|^2_{B',0,8}\\
&\leq&C\|\psi\|_{B,0,4},
\end{eqnarray*}
$$\||\psi|^3\|_{B,0,p}\leq C\|\psi\|_{B,0,4},$$
$$\||\psi|^5\|_{B,0,p}\leq C\|\psi\|_{B,0,4}.$$
Substituting these into (\ref{2.18}) on $B'$ with $p=\frac 4 3$
and using (\ref{2.13}), (\ref{2.17}) again, we conclude that
$$\|\eta\psi\|_{B',2,\frac 4 3}\leq C\|\psi\|_{B,0,4}.$$
This implies that
$$\|\psi\|_{B'',2,\frac 4 3}\leq C\|\psi\|_{B,0,4}$$
for any $B''\subset \subset B'$. By Sobolev,
$$\|\psi\|_{B'',1, 4}\leq C\|\psi\|_{B,0,4},$$
and consequently,  $\|\psi\|_{L^\infty(B'')}\leq C\|\psi\|_{B,0,4}.$
\par
Choose $p=2$ in (\ref{2.18}) and use the above estimates, we have
$$\|\psi\|_{B',2, 2}\leq C\|\psi\|_{B,0,4}.$$
This yields
$$\|\psi\|_{B',1, p}\leq C\|\psi\|_{B,0,4}.$$
We can then obtain all the desired estimates by the standard
bootstrap method, for example, using $\|\psi\|_{B',1, p}\leq
C\|\psi\|_{B,0,4}$ in (\ref{2.18}) we have $\|\psi\|_{B'',2, p}\leq
C\|\psi\|_{B,0,4}$ for any $B''\subset \subset B'$ and $p>1$. \hfill
Q.E.D. \vskip12pt From the above estimates, we have the following
(see e.g. \cite{MS}): \vskip12pt \noindent {\bf Corollary 2.3.} {\it
Let $M$ be a compact Riemannian surface without boundary, with a
fixed spin structure. If $\psi\in W^{1,p}$ is a solution of
\begin{equation}
\label{2.19}
\D\psi=F
\end{equation}
with $F\in L^p$, then
\begin{equation}
\label{2.20}
\|\psi\|_{M,1,p}\leq C(\|F\|_{M,0,p}+\|\psi\|_{M,0,p}),
\end{equation}
where $C=C(M,p)>0$ is a constant.  } \vskip12pt \noindent {\it
Proof.}  Assume that $\{U_\alpha, \varphi_\alpha\}$ is a finite
covering of $M$ by charts. Let $\{g_\alpha\}$ be a partition of unit
subordinate to this covering. Denote $\psi_\alpha:=g_\alpha \psi$,
and $F_\alpha:=\D\psi_\alpha=g_\alpha F+\nabla g_\alpha\cdot\psi$,
then in each chart $\{U_\alpha, \varphi_\alpha\}$, using (\ref{2.7})
we have
$$\|\psi_\alpha\|_{U_\alpha,1,p}\leq
C_\alpha(\|F_\alpha\|_{M,0,p}+\|\psi_\alpha\|_{M,0,p}),$$ noting
that $\psi=\sum\limits_\alpha\psi_\alpha$, $\|\nabla g_\alpha\|\leq
C$ and $g_\alpha<1$, we obtain (\ref{2.20}). \hfill Q.E.D.

\vskip36pt
\section{Removable singularity theorem}
\addtocounter{equation}{-20}

For a given smooth solution $\psi$ of (\ref{1.1}) on the sphere
$\Bbb{S}^2$, one can create a solution $\widetilde{\psi}$ on the
Euclidean plane $\R^2$ through the stereographic projection from the
north pole $N$, by virtue of the conformal invariance of the
equation. Conversely, given a solution $\psi$ on $\R^2$, through the
stereographic projection, we only have a solution $\psi$ on
$\Bbb{S}^2\setminus \{N\}$, which then leads to the question of
removable singularities. In this section, we will prove the
following

\vskip12pt \noindent {\bf Theorem 3.1 (Removable singularity
theorem).} {\it Let $\psi$ be a solution of (\ref{1.1}) which is
smooth on $B\setminus \{0\}$. If
\begin{equation}
\label{3.1}
\int\limits_B|\psi|^4 <\infty,
\end{equation}
then $\psi$ extends to a smooth solution of (\ref{1.1}) on the whole $B$.}
\vskip12pt
\noindent
{\it Proof.}  Since (\ref{1.1}) is conformally invariant, by a rescaling transformation, we may assume that
$$\int\limits_B|\psi|^4 <\varepsilon,$$
where $\varepsilon>0$ is a small constant whose appropriate value will
be determined later. For any given small $\delta>0$, we choose a cut-off
function $\eta_\delta\in C^\infty_0(B_{2\delta}) $  such that $0\leq \eta_\delta\leq 1$,
$$\eta_\delta=\left\{
                \begin{array}{l}
  1\quad {\rm in} \quad B_\delta, \\
  0\quad {\rm in}\quad B\setminus B_{2\delta},
                \end{array}
                \right.
$$
and $|\nabla \eta_\delta|\leq  C/\delta$. Then
\begin{eqnarray*}
\D[(1-\eta_\delta)\psi]&=&(1-\eta_\delta)\D\psi-\nabla\eta_\delta\cdot\psi\\
&=&(1-\eta_\delta)H_{jkl}\la\psi^j,\psi^k\ra\psi^l-\nabla\eta_\delta\cdot\psi.
\end{eqnarray*}
By Lemma 2.2, we have
\begin{eqnarray}
\qquad \|(1-\eta_\delta)\psi\|_{B,1,\frac 4 3}&\leq
&C\|(1-\eta_\delta)H_{jkl}\la\psi^j,\psi^k\ra\psi^l-\nabla\eta_\delta\cdot\psi\|_{B,0,\frac
4 3}
+C\|\psi\|_{\partial B,1,\frac 4 3}\\
&\leq& C (h_0\|\psi\|^3_{B,0,4}+\|\nabla\eta_\delta\cdot\psi\|_{B,0,\frac 4 3}+\|\psi\|_{\partial B,1,\frac 4 3}).
\end{eqnarray}
By the Sobolev embedding theorem, we have
\begin{equation}
\label{3.3} \|(1-\eta_\delta)\psi\|_{B,0,4}\leq
C(h_0\|\psi\|^3_{B,0,4}+\|\nabla\eta_\delta\cdot\psi\|_{B,0,\frac
4 3}+\|\psi\|_{\partial B,1,\frac 4 3}).
\end{equation}
We note that as $\delta \to 0$,
\begin{eqnarray*}
\|\nabla\eta_\delta\cdot\psi\|_{B,0,\frac 4
3}&=&(\int\limits_{B_{2\delta}\setminus
B_\delta}|\nabla\eta_\delta|^{\frac 4 3}|\psi|^{\frac 4 3})^{\frac
3 4}\\
&\leq&\frac C \delta (\int\limits_{B_{2\delta}}|\psi|^{\frac 4
3})^{\frac 3 4}\\
&\leq&C (\int\limits_{B_{2\delta}}|\psi|^4 )^{\frac 1 4}\to 0,
\end{eqnarray*}
therefore, letting $\delta \to 0$ in (\ref{3.3}) we obtain
$$\|\psi\|_{B,0,4}\leq Ch_0\|\psi\|^2_{B,0,4}\|\psi\|_{B,0,4}+C\|\psi\|_{\partial B,1,\frac 4
3}.$$
We choose $\varepsilon>0$ so small  that
$Ch_0\sqrt{\varepsilon}<1/2$, then
\begin{eqnarray*}
\|\psi\|_{B,0,4}&\leq&C\|\psi\|_{\partial B,1,\frac 4 3}\\
&\leq&C(\int\limits_{\partial B}|\nabla\psi|^{\frac 4 3})^{\frac 3
4}+C(\int\limits_{\partial B}|\psi|^4)^{\frac 1 4}.
\end{eqnarray*}
By a rescaling argument, we have for any $r\in (0,1]$,
\begin{eqnarray*}
(\int\limits_{B_r}|\psi|^4)^{\frac 1 4}
&\leq&C(r\int\limits_{\partial B_r}|\nabla\psi|^{\frac 4
3})^{\frac 3 4}+C(r\int\limits_{\partial B_r}|\psi|^4)^{\frac 1
4}\\
&\leq&C[(r\int\limits_{\partial B_r}|\nabla\psi|^{\frac 4
3})^{\frac 1 4}+(r\int\limits_{\partial B_r}|\psi|^4)^{\frac 1
4}],
\end{eqnarray*}
that is,
\begin{equation}
\label{3.4} \int\limits_{B_r}|\psi|^4\leq Cr\int\limits_{\partial
B_r}|\nabla\psi|^{\frac 4 3}+Cr\int\limits_{\partial B_r}|\psi|^4.
\end{equation}
Denote
$$\bar{\psi}:=\frac{1}{2\pi}\int\limits_{\partial B}\psi,$$
then on $B\setminus\{0\}$:
$$\D(\psi-\bar{\psi})=H_{jkl}\la\psi^j,\psi^k\ra\psi^l=H_{jkl}\la\psi^j,\psi^k\ra(\psi^l-\bar{\psi}^l)+H_{jkl}\la\psi^j,\psi^k\ra\bar{\psi}^l.$$
From Lemma 2.2, we have
\begin{eqnarray*}
\|\psi-\bar{\psi}\|_{B,1,\frac 4
3}&\leq&C(\|H_{jkl}\la\psi^j,\psi^k\ra(\psi^l-\bar{\psi}^l)\|_{B,0,\frac
4
3}+\|H_{jkl}\la\psi^j,\psi^k\ra\bar{\psi}^l\|_{B,0,\frac 4 3}\\
&&+\|\psi-\bar{\psi}\|_{\partial B,1,\frac 4 3}),
\end{eqnarray*}
using the Poincare's inequality, we obtain
\begin{eqnarray*}
\|\psi-\bar{\psi}\|_{B,1,\frac 4 3}&\leq &
C(\|H_{jkl}\la\psi^j,\psi^k\ra(\psi^l-\bar{\psi}^l)\|_{B,0,\frac 4
3}+\|H_{jkl}\la\psi^j,\psi^k\ra\bar{\psi}^l\|_{B,0,\frac 4
3}\\
&&+\|\nabla(\psi-\bar{\psi})\|_{\partial B,0,\frac 4 3})\\
&\leq&
Ch_0\|\psi\|^2_{B,0,4}\|\psi-\bar{\psi}\|_{B,0,4}+Ch_0\|\psi\|^2_{B,0,4}\|\bar{\psi}\|_{B,0,4}
\\
&&+C\|\nabla\psi\|_{\partial B,0,\frac 4 3}\\
&\leq&Ch_0\|\psi\|^2_{B,0,4}\|\psi-\bar{\psi}\|_{B,1,\frac 4
3}+Ch_0\|\psi\|^2_{B,0,4}\|\psi\|_{\partial B,0,1}\\
&& +C\|\nabla\psi\|_{\partial B,0,\frac 4 3},
\end{eqnarray*}
using the smallness of $\|\psi\|_{B,0,4}$ again, we have
$$\|\psi-\bar{\psi}\|_{B,1,\frac 4 3}\leq Ch_0\|\psi\|^2_{B,0,4}\|\psi\|_{\partial B,0,4}+ C\|\nabla\psi\|
_{\partial B,0,\frac 4 3},$$
 therefore,
 $$ \|\nabla\psi\|_{B,0,\frac 4 3}\leq Ch_0\|\psi\|^2_{B,0,4}\|\psi\|_{\partial B,0,4}+ C\|\nabla\psi\|
_{\partial B,0,\frac 4 3},  $$ that is,
$$(\int\limits_B|\nabla\psi|^{\frac 4 3})^{\frac 3 4}\leq
Ch_0(\int\limits_B|\psi|^4)^{\frac 1 2}(\int\limits_{\partial
B}|\psi|^4)^{\frac 1 4}+C(\int\limits_{\partial
B}|\nabla\psi|^{\frac 4 3})^{\frac 3 4}.$$
It is then easy to see
that
\begin{eqnarray*}
\int\limits_B|\nabla\psi|^{\frac 4 3}&\leq&Ch_0^{\frac 4
3}(\int\limits_B|\psi|^4)^{\frac 2 3}(\int\limits_{\partial
B}|\psi|^4)^{\frac 1 3}+C\int\limits_{\partial
B}|\nabla\psi|^{\frac 4 3}\\
&\leq& \sigma\int\limits_B|\psi|^4+\frac
{C}{\sigma}\int\limits_{\partial B}|\psi|^4+C\int\limits_{\partial
B}|\nabla\psi|^{\frac 4 3},
\end{eqnarray*}
where $\sigma>0$ is small constant, and $C>0$ is constant depending
only on $h_0$. By rescaling again, we have
\begin{equation}
\label{3.5} \int\limits_{B_r}|\nabla\psi|^{\frac 4 3}\leq
\sigma\int\limits_{B_r}|\psi|^4+\frac
{C}{\sigma}r\int\limits_{\partial
B_r}|\psi|^4+Cr\int\limits_{\partial B_r}|\nabla\psi|^{\frac 4 3}
\end{equation}
Combining (\ref{3.4}) and (\ref{3.5}), we have
\begin{equation}
\label{3.6}
\int\limits_{B_r}|\psi|^4+\int\limits_{B_r}|\nabla\psi|^{\frac 4
3}\leq Cr(\int\limits_{\partial B_r}|\psi|^4+\int\limits_{\partial
B_r}|\nabla\psi|^{\frac 4 3}).
\end{equation}
Denote $F(r):=\int\limits_{B_r}|\psi|^4+|\nabla\psi|^{\frac 4 3}$,
then
$$F(r)\leq CrF'(r),$$
which implies that
\begin{equation}
\label{3.7} F(r)\leq F(1)r^{\frac 1 C}.
\end{equation}
From this, it follows that $\psi\in W^{1,p}$ for some $p>4/3$, and
then,
by the standard bootstrap method, one can conclude the smoothness
of $\psi$. \hfill Q.E.D.

\vskip12pt \noindent{\bf Remark.} When $\psi$ is the
 spinor representing a surface $M$ in $\R^3$ with mean
curvature $H$, and $z=x+iy$ is the parameterization of $M$, then the metric of $M$ is
$$ds^2=|\psi|^4dzd\bar{z},$$
and the condition (\ref{3.1}) means that $(M,ds^2)$ has finite area.

\vskip36pt
\section{Energy identity}
\addtocounter{equation}{-8}

Let $M$ be a compact Riemann surface with a fixed spin structure. Given a sequence $\{\psi_m\}$ of solutions of (\ref{1.1}) on $M$, if
we assume it is uniformly bounded in $L^p(p>4)$, then the standard
bootstrap method implies that $\{\psi_m\}$ is uniformly bounded in
$C^r$ $(r\in \Z_+)$. However, in the case of the $L^4$-norm, examples
show that this compactness is no longer true. If $\{\psi_m\}$
converges to $\psi$ weakly in $L^4$, then in the limit we may encounter
bubbling phenomenon; namely, by a rescaling argument and the previous
removable singularity theorem, we may get some solutions on
$\Bbb{S}^2$, and this causes an energy loss. Comparing to the well-known
case of harmonic maps, one naturally asks whether the blow up set is
finite and the energy identity holds. In view of the Weierstrass
representation, this corresponds to the question of the
convergence of surfaces with a uniform area bound.
\par
 We will need the following lower
bound for the energy of the bubbles:

\vskip12pt \noindent {\bf Lemma 4.1.} {\it There exists a constant
$A>0$ such that for any nontrivial solution $\psi$ of (\ref{1.1}) on
$\Bbb{S}^2$, we have
\begin{equation}
\label{4.1}  \int\limits_{\Bbb{S}^2}|\psi|^4\geq A.
\end{equation}
}
 \noindent {\it Proof.} Firstly,  for any
solution $\psi$ on $\Bbb{S}^2$,
\begin{equation}
\label{4.2} \|\psi\|_{\frac 4 3}\leq C\|\D\psi\|_{\frac 4 3}.
\end{equation}
Otherwise, for any $k\in \Z_+$, there is a $\psi_k$ which solves
(\ref{1.1}), but
$$\|\psi_k\|_{\frac 4 3}>k\|\D\psi_k\|_{\frac 4 3}.$$
Denote $\eta_k:=\psi_k/\|\psi_k\|_{\frac 4 3 }$, then
\begin{equation}
\label{4.3} \|\D\eta_k\|_{\frac 4 3}<1/k, \quad \|\eta_k\|_{\frac
4 3}=1.
\end{equation}
Using Corollary 2.3, we have
$$\|\eta_k\|_{1,\frac 4 3}\leq C(\|\D\eta_k\|_{\frac 4
3}+\|\eta_k\|_{\frac 4 3})\leq C,$$ which implies that there exists
some $\eta_0$ such that $\eta_k$ converges to $\eta_0$ weakly in
$W^{1,\frac 4 3}$. By Sobolev, $\eta_k\to \eta_0$ in $L^{\frac 4
3}$, so $\|\eta_0\|_{\frac 4
3}=\lim\limits_{k\to\infty}\|\eta_k\|_{\frac 4 3}=1.$ But $\D\eta_k$
converges to $\D\eta_0$ weakly in $L^{\frac 4 3}$, and from
(\ref{4.3}), it is easy to see that $\D\eta_0=0$, hence
$\eta_0\equiv 0$ since there is no nontrivial harmonic spinor on
$\Bbb{S}^2$. This contradicts  $\|\eta_0\|_{\frac 4 3}=1.$

Now from (\ref{4.2}) and (\ref{2.20}), we have
\begin{eqnarray*}
\|\psi\|_{1,\frac 4 3}&\leq& C(\|\D\psi\|_{\frac 4
3}+\|\psi\|_{\frac 4 3})\\
&\leq &C\|\D\psi\|_{\frac 4 3}\\
&\leq & Ch_0\|\psi\|^3_4,
\end{eqnarray*}
therefore, $$\|\psi\|_4\leq Ch_0\|\psi\|^3_4,$$ and if $\|\psi\|_4$ is
so
small  that $Ch_0\|\psi\|^2_4<1$, then we have
$\psi\equiv 0$. Equivalently, we can find a constant $A>0$ such that
for any nontrivial solution $\psi$ of (\ref{1.1}) on $\Bbb{S}^2$,
the energy $\int\limits_M|\psi|^4$ is bounded below by $A$. \hfill
Q.E.D.

\par
Let
$$S:=\cap_{r>0}\{x\in M|\liminf\limits_{m\to
+\infty}\int_{B(x,r)}|\psi_m|^4\geq\varepsilon\}$$  be the blow up
set of $\{\psi_m\}$, where $\varepsilon$ is as in Theorem 2.1.

\vskip12pt \noindent {\bf Theorem 4.2 (Energy Identity).} {\it Let
$M$ be a compact Riemann surface with fixed spin structure, and
suppose that $\{\psi_m\}$ is a sequence of smooth solutions of
(\ref{1.1}) on $M$ satisfying
\begin{equation}
\label{4.4} E(\psi_m):=\int\limits_M|\psi_m|^4\leq \Lambda<+\infty.
\end{equation}
If $\{\psi_m\}$ converges to $\psi$ weakly in $L^4(M)$ (but not
strongly,) then the (non-empty) blow up set $S$ must be finite:
$$S=\{p_1,p_2,\cdots,p_K\}.$$
Furthermore, there exists a constant $c_0>0$ depending only on $M$
such that if
\begin{equation}
\label{4.5.1} \sup\limits_{M,i,j,k,l}|H^i_{jkl}|\sqrt{\Lambda}<c_0,
\end{equation}
 then the energy identity for
$\{\psi_m\}$ holds, namely, for each blow up point $p_k$
$(k=1,2,\cdots,K)$, there exist a finite number of solutions
$\{\xi^a_k\}_{a=1,2,\cdots,A_k}$ of (\ref{1.1}) on $\Bbb{S}^2$ such
that
\begin{equation}
\label{4.5} \lim\limits_{n\to
+\infty}E(\psi_m)=E(\psi)+\sum\limits_{k=1}^K\sum\limits_{a=1}^{A_k}E(\xi^a_k).
\end{equation}
 }
 \noindent {\it Proof.} Since the removable singularity
theorem, the small energy regularity theorem, and Lemma 4.1 provide
key ingredients for establishing the energy identity, the theorem
can then be proved by an argument as the proof of Theorem 3.6 in
\cite{CJLW2}, see also \cite{DT}. Here we only give a sketch of
proof.

First, the condition $E(\psi_m)\leq\Lambda<+\infty$ and Theorem 2.1
imply that the blow up set $S$ must be finite. We choose small disks
$B_{\delta_k}$ for each $p_k$ such that $B_{\delta_k}\cap
B_{\delta_j}=\phi$ for $k\not= j,$ $ k,j=1,2,\cdots,K.$ Furthermore,
by Theorem 2.1, $\{\psi_m\}$ strongly converges to $\psi$ in $L^4$
on $M\setminus\cup_{k=1}^K B_{\delta_k},$ (\ref{4.5}) is then
equivalent to
\begin{equation}
\label{c7} \Sigma^K_{k=1} \lim \limits_{\delta_k\to 0}\lim
\limits_{n\to\infty}E(\psi_m;B_{\delta_k})=\Sigma^K_{k=1}
\Sigma_{a=1}^{A_k}E(\xi^a_k).
\end{equation}
It suffices to prove that for each blow-up point $p$, we have
\begin{equation}
\label{c9} \lim \limits_{\delta\to 0}\lim
\limits_{m\to\infty}E(\psi_m;B_{\delta})= \Sigma_{a=1}^{A}E(\xi^l).
\end{equation}
By virtue of the conformal invariance of the equation (\ref{1.1}),
and the locality of the problem, we may assume that each disk
$B_\delta$ is equipped with the Euclidean metric. For each $\psi_m$,
we choose $\lambda_m$ and $x_m\in B_\delta$ such that $\lambda_m \to
0$, $x_m \to p$ and

$$E(\psi_m;
B_{\lambda_m}(x_m))=\max\limits_{x\in B_\delta(p)} E(\psi_m;
B_{\lambda_m}(x))=\frac{\varepsilon}{2}.$$

\noindent Rescaling by
$$ \widetilde{\psi}_m(x):=\lambda_m^{-\frac 1 2 }
\psi_m(x_m+\lambda_mx),$$ then
$$E(\widetilde{\psi}_m;B)=E(\psi_m;B_{\lambda_m}(x_m))
=\frac{\varepsilon}{2}<\varepsilon,$$
$$E(\widetilde{\psi}_m;B_R)=E(\psi_m;B_{\lambda_m R}(x_m))\leq
\Lambda.$$
By Theorem 2.1, we have a subsequence of $\{\psi_m\}$
which strongly converges to some $\widetilde{\psi}$ in $L^4(B_R)$
for any $R\geq 1$. We thus obtain a nonconstant solution
$\widetilde{\psi}$ of (\ref{1.1}) on $\RR^2$, and hence a
nonconstant solution of (\ref{1.1}) $\xi^1$ on
$\Bbb{S}^2\setminus\{N\}$ with bounded energy. Theorem 3.1 then
gives us a nonconstant solution of (\ref{1.1}) on the whole
$\Bbb{S}^2,$ and we obtain the first bubble $\xi^1$ at the blow up point
$p$.

\par
Next, denote $$A(\delta, R,m):=\{x\in \RR^2|\lambda_m R\leq
|x-x_m|\leq\delta\},$$ then (\ref{c9}) is equivalent to

\begin{equation}
\label{c11} \lim \limits_{R\to\infty}\lim \limits_{\delta\to 0}\lim
\limits_{m\to \infty}E(\psi_m;A(\delta,R,m))=\Sigma_{a=2}^A
E(\xi^a).
\end{equation}

For a fixed blow-up point $p$, the number of bubbles $\xi$ must be
finite; this follows easily from Lemma 4.1.
 We only consider (\ref{c11}) in the case of exactly one bubble at the
blow up point $p$, because the case of at least two bubbles can be reduced
to this case. Then, (\ref{c11}) is just

\begin{equation}
\label{c13} \lim \limits_{R\to\infty}\lim \limits_{\delta\to 0}\lim
\limits_{m\to \infty}E(\psi_m;A(\delta,R,m))=0.
\end{equation}
To prove this, we consider a conformal transformation $f: \RR\times
\Bbb{S}^1\to \RR^2, \quad f(t,\theta)=(e^{-t},\theta)$, where
$\RR\times \Bbb{S}^1$ is given the metric $g=dt^2+d\theta^2$. for
the pull-back $\Psi_m:=f^*\psi_m,$ then $E(\Psi_m)\leq \Lambda$. Set
$T_0:=|{\rm log}\delta|,$ $T_m:=|{\rm log}\lambda_m R|.$

Using Theorem 3.1, through an argument by contradiction
  (c.f. p.82 in \cite{CJLW2} for more details), one can prove that there is a
$K>0$ such that if $m\geq K$, then
\begin{equation}
\label{c16} \int_{[t,t+1]\times\Bbb{S}^1}|\Psi_m|^4<\varepsilon,
\quad \forall t\in [T_0,T_m-1].
\end{equation}

\par
Choose a cut-off function $\eta$ on $B(x_m,2\delta)$ as follows:
$$\eta\in C^\infty_0(B_{2\delta}\setminus
B_{\lambda_mR/2});\qquad \eta\equiv 1 \quad {\rm in}\quad
B_{\delta}\setminus B_{\lambda_m R}$$
$$|\nabla\eta|\leq C/\delta \quad {\rm in}\quad B_{2\delta}\setminus
B_\delta;\qquad |\nabla\eta|\leq C/\lambda_mR \quad {\rm in}\quad
B_{\lambda_mR}\setminus B_{\lambda_mR/2},$$ where we denote
$B_\delta:=B(x_m,\delta)$ etc. for simplicity.
 Then from Lemma 2.2 we have
\begin{eqnarray*}
\|\eta\psi_m\|_{L^4}&\leq&C\|\eta \D
\psi_m+\nabla\eta\cdot\psi_m\|_{L^{\frac 4 3}} \\
&\leq&C\|h_0|\eta||\psi_m|^3+|\nabla\eta||\psi_m|\|_{L^{\frac 4 3}}\\
&\leq& Ch_0\|\psi_m\|^2_{L^4}\|\eta\psi_m\|_{L^4}
+C[\int_{A(2\delta,R/2,m)}(|\nabla\eta||\psi_m|)^{\frac 4 3}]^{\frac 3 4}\\
&\leq& Ch_0\sqrt{\Lambda}\|\eta\psi_m\|_{L^4}
+C[\int_{A(2\delta,R/2,m)}(|\nabla\eta||\psi_m|)^{\frac 4 3}]^{\frac
3 4}.
\end{eqnarray*}
Clearly, there exists a constant $c_0>0$ such that when (\ref{4.5.1}) is satisfied, we have
$Ch_0\sqrt{\Lambda}<1$, from the above estimate, we then have
\begin{eqnarray*}
\|\eta\psi_m\|_{L^4} &\leq& C[\int_{B_{2\delta}\setminus
B_\delta}(|\nabla\eta| |\psi_m|)^{\frac 4 3}]^{\frac 3
4}+C[\int_{B_{\lambda_mR}\setminus
B_{\lambda_mR/2}}(|\nabla\eta||\psi_m|)^{\frac 4 3}]^{\frac 3 4}.
\end{eqnarray*}
Therefore,
\begin{eqnarray*}
\|\psi_m\|_{L^4(A(\delta,R,m))}&\leq& C [\int_{B_{2\delta}\setminus
B_\delta}
|\psi_m|^4]^{\frac 1 4}+C[\int_{B_{\lambda_m R}\setminus B_{\lambda_mR/2}}|\psi_m|^4]^{\frac 1 4}\\
&\leq&C\varepsilon^{\frac 1
4}+C\varepsilon^{\frac 1 4},
\end{eqnarray*}
where in the last step, we used (\ref{c16}).  This proves
(\ref{c13}). \hfill Q.E.D.

\vskip36pt
\section{Related types of Dirac equations}
\addtocounter{equation}{-11}

Let $M$ be a compact Riemann surface with fixed spin structure. For any local orthonormal basis $\{e_\alpha\}_{\alpha=1,2}$, one can
define the so-called chirality operator
$$\Gamma:=i\hskip0.8mm e_1\cdot e_2\cdot .$$
This definition is independent of the choice of
$\{e_\alpha\}_{\alpha=1,2}$, therefore $\Gamma$  is globally defined
on $M$. Define
$$ \Gamma_+:=\frac 1 2 (Id+\Gamma),  \quad
           \Gamma_-:=\frac 1 2 (Id-\Gamma).$$
Let $U=U(\psi)$, $V=V(\psi)$ be complex functions. We consider the
following Dirac equation:
\begin{equation}
\label{5.1} \D\psi=[U(\psi)\Gamma_++V(\psi)\Gamma_-]\psi.
\end{equation}
Clearly, (\ref{1.1}) corresponds to the case $U=V=-H|\psi|^2$.
Comparing to (\ref{1.1}), the Dirac equation for surfaces immersed
into some three-dimensional Lie group $N$ takes a special form of
(\ref{5.1}). For example (c.f.  section 2.3 in \cite{Ta2})
\begin{equation}
\label{5.2}
 N=SU(2):\qquad   U=\bar{V}=-(H-i)|\psi|^2;
\end{equation}
\begin{equation}
\label{5.3}
 N=Nil:\qquad   U=V=-H|\psi|^2-\frac i 2 (|\psi_1|^2-|\psi_2|^2);
\end{equation}
\begin{eqnarray}
\label{5.4}
 N=\widetilde{SL_2}:\quad   &&U=-H|\psi|^2-i(\frac 3 2 |\psi_2|^2-|\psi_1|^2), \nonumber\\
 &&V=-H|\psi|^2-i( |\psi_2|^2-\frac 3 2|\psi_1|^2).
\end{eqnarray}
Using the same methods as in the previous sections, one can conclude
similar results for the above types of Dirac equations.

\vskip12pt\noindent {\bf Corollary 5.1.} {\it  There exists a
small constant $\varepsilon>0$ such that for any smooth solution
$\psi$ of (\ref{5.1}), with $U,V$ satisfying one of
(\ref{5.2}),(\ref{5.3}),(\ref{5.4}), and
\begin{equation}
\label{5.5} E(\psi;B):=\int\limits_B|\psi|^4<\varepsilon,
\end{equation}
we have
\begin{equation}
\label{5.6}
\|\psi\|_{B',k,p}\leq C\|\psi\|_{B,0,4},
\end{equation}
$\forall B'\subset \subset B,$ $1<p$ and $k\in\Z_+,$ where
$C=C(B',k,p)>0$ is constant, and $\|\cdot\|_{B,k,p}$ denotes the
norm in $W^{k,p}(B,\S^n)$. } \vskip12pt\noindent {\bf Corollary
5.2.} {\it Let $\psi$ be a solution of (\ref{5.1}),with $U,V$
satisfying one of (\ref{5.2}),(\ref{5.3}),(\ref{5.4}), which is
smooth on $B\setminus \{0\}$. If
\begin{equation}
\label{5.7}
\int\limits_B|\psi|^4 <\infty,
\end{equation}
then $\psi$ extends to a smooth solution  on the whole $B$.
}

\vskip12pt\noindent {\bf Corollary 5.3.} {\it  Let $M$ be a compact
Riemann surface with fixed spin structure, and suppose that
$\{\psi_m\}$ is a sequence of smooth solutions of (\ref{5.1}), with
$U,V$ satisfying one of (\ref{5.2}),(\ref{5.3}),(\ref{5.4})
respectively, on $M$ and
\begin{equation}
\label{5.8} E(\psi_m):=\int\limits_M|\psi_m|^4\leq \Lambda<+\infty.
\end{equation}
If $\{\psi_m\}$ converges to $\psi$ weakly in $L^4(M)$ (but not
strongly), then the (non-empty) blow up set $S$ must be finite:
$$S=\{p_1,p_2,\cdots,p_K\}.$$
Furthermore, there exists a constant $c_0>0$ depending only on $M$
such that if
\begin{equation}
\label{5.9} (\sup\limits_M|H|+\alpha)\sqrt{\Lambda}<c_0,
\end{equation}
with $\alpha=1, \frac 1 2, \frac 3 2$ respectively,
 then the energy identity for
$\{\psi_m\}$ holds, namely, for each blow up point $p_k$
$(k=1,2,\cdots,K)$, there exist a finite number of solutions
$\{\xi^a_k\}_{a=1,2,\cdots,A_k}$ of (\ref{5.1}) on $\S^2$ such that
\begin{equation}
\label{5.10} \lim\limits_{m\to
+\infty}E(\psi_m)=E(\psi)+\sum\limits_{k=1}^K\sum\limits_{a=1}^{A_k}E(\xi^a_k).
\end{equation}
}

\end{document}